\documentclass[12pt,a4paper]{article}

\RequirePackage{geometry}
\usepackage{amssymb,amsmath}
\usepackage{graphicx}
\usepackage{graphics}
\usepackage{hyperref}
\usepackage{subfigure}
\usepackage{mathrsfs}
\usepackage{eucal}
\usepackage{amsmath}
\usepackage{amssymb}
\usepackage{amsfonts}
\usepackage{multirow}
\usepackage{lscape}

\usepackage{natbib}

\newcommand{\fref}[1]{Figure (\ref{#1})}

\newcommand{\KK}{\mathbf{K}}
\newcommand{\bm}{\mathbf{M}}

\newcommand{\DD}{\mathbf{D_b}}

\newcommand{\bveps}{\boldsymbol{\varepsilon}}

\newcommand{\BB}{\mathbf{B}}

\newcommand{\rmd}{\mathrm{d}}

\begin{document}

\begin{center}
\large{Supersonic flutter analysis of thin cracked functionally graded material plates}
\end{center}

\begin{center}
S~Natarajan$^{1,a}$, M~Ganapathi$^2$, S~Bordas$^3$
\end{center}

\begin{center}\small{
$^{1,a}$Institute of Mechanics and Advanced Materials, Theoretical and Computational Mechanics, Cardiff University, U.K. Email: sundararajan.natarajan@gmail.com. \\
$^2$Head, Stress \& DTA, Bombardier Aerospace India Center, IES-Aerospace, Mahindra Satyam Computers Services Ltd., Bangalore, India \\
$^3$Professor, Institute of Mechanics and Advanced Materials, Theoretical and Computational Mechanics, Cardiff University, U.K.}
\end{center}

\begin{abstract}
In this paper, the flutter behaviour of simply supported square functionally graded material plates immersed in a supersonic flow is studied. An enriched 4-noded quadrilateral element based on field consistency approach is used for this study and the crack is modelled independent of the underlying mesh. The material properties are assumed to be temperature dependent and graded only in the thickness direction. The effective material properties are estimated using the rule of mixtures. The formulation is based on the first order shear deformation theory and the shear correction factors are evaluated employing the energy equivalence principle. The influence of the crack length, the crack orientation, the flow angle and the gradient index on the aerodynamic pressure and the frequency are numerically studied. The results obtained here reveal that the critical frequency and the critical pressure decreases with increase in crack length and it is minimum when the crack is aligned to the flow angle.
\end{abstract}

\begin{footnotesize}
\textbf{Keywords}: Mindlin plate theory, flutter, partition of unity methods, extended finite element method, shear flexible element, functionally graded material.
\end{footnotesize}

\section{Introduction}
The emergence of functionally graded materials (FGMs) that combine the best properties of its constituent materials (for example, ceramics and metals) is considered to be an alternative for certain class of aerospace structures exposed to high temperature environment. FGMs are characterized by a smooth transition from one material to another, thus circumventing high inter-laminar shear stresses and de-lamination problem that persists in laminated composites. The introduction of FGMs has attracted researchers to investigate the structural behaviour of such structures. With the increased use of these materials, it necessitates to study the dynamic characteristics of the structures made up of FGMs. 

The investigation of dynamic behaviour of FGM structures is fairly dealt in the literature. Some of the important contributions are discussed here. He \textit{et al.,}~\citep{He2001} presented finite element formulation based on thin plate theory for the vibration control of FGM plate with integrated piezoelectric sensors and actuators under mechanical load, whereas Liew \textit{et al.,}~\citep{Liew1994} have analyzed the active vibration control of plate subjected to a thermal gradient using shear deformation theory. The parametric resonance of FGM plates is dealt by Ng \textit{et al.,}~\citep{Ng2000} based on Hamilton's principle and the assumed mode technique. Yang and Shen~\citep{Yang2001,Yang2002} have analyzed the dynamic response of thin FGM plates subjected to impulsive loads using Galerkin Procedure coupled with modal superposition methods, whereas, by neglected the heat conduction effect, such plates and panels in thermal environments have been examined based on shear deformation with temperature dependent material properties. The static deformation and vibration of FGM plates based on higher-order shear deformation theory is studied by Qian \textit{et al.,}~\citep{Qian2004a} using meshless local Petrov-Galerkin method (MLPG). Matsunaga~\citep{Matsunaga2008} presented analytical solutions for simply supported rectangular FGM plates based on second-order shear deformation plates, whereas, three dimensional solutions are proposed in~\citep{Vel2002,Vel2004} for vibrations of simply supported rectangular FGM plates. Reddy~\citep{Reddy2000} presented finite element solutions for the dynamic analysis of a FGM plate and Ferreira \textit{et al.,}~\citep{Ferreira2006} performed dynamic analysis of FGM plate based on higher order shear and normal deformable plate theory using MLPG. Since FGMS are seen as potential candidates for aircraft structural applications, it is important to understand the dynamic characteristics of structures made up of such materials when exposed to air flow. Some of the recent contributions are discussed here. Prakash and Ganapathi~\citep{Prakash2006} studied the linear flutter characteristics of FGM panels exposed to supersonic flow. Haddadpour \textit{et al.,}~\citep{Haddadpour2007} and Sohn and Kim~\citep{Sohn2008,Sohn2009} investigated the nonlinear aspects of flutter characteristics using the finite element method.

FGM plates or in general plate structures may develop flaws during manufacturing or while in service. The static and dynamic fracture mechanics study of FGMs have been studied in~\citep{Dolbow2002,Nazari2011}. Dolbow and Gosz~\citep{Dolbow2002} employed extended finite element method (XFEM) to compute mixed mode stress intensity factors. Dynamic analysis of FGM beams with edge cracks has been attempted by many researchers~\citep{Yang2008,Kitipornchai2009,Yan2011,Yan2011a}. However, such analysis of plates is scarce in the literature~\citep{Huang2011,natarajanbaiz2011,Natarajan2011}. Huang \textit{et al.,}~\citep{Huang2011} have analyzed the vibration of side-cracked FGM thick plate analytically by employing Ritz procedure, whereas, Natarajan \textit{et al.,}~\citep{natarajanbaiz2011,Natarajan2011} examined the FGM plate with through center crack using the XFEM. It can be seen from the available literature that the work on flutter characteristics of FGM plates with cracks, to the author's knowledge is not available in the literature. Earlier studies on flutter characteristics of cracked isotropic and composite panels employed finite element procedure~\citep{Chen1985,Pidaparti1998,Shiau1992,Strganac1996}. Although these numerical studies give insight into understanding the flutter behaviour, the method requires the mesh to conform to the geometry. This inherently limits the analyses to fixed calculation parameters. With the introduction of the extended finite element method (XFEM)~\citep{Belytschko1999}, the requirement of mesh conforming to the geometry is alleviated and it is now possible to model and/or track irregular geometries as they evolve. 

In this paper, we apply the XFEM~\citep{Belytschko1999} to model the crack independent of the underlying mesh and then to study the flutter behaviour of FGM plates immersed in a supersonic flow based on first order shear deformation theory. The crack kinematics is captured by adding suitable functions to the finite element basis. These additional functions are called the `enrichment functions' and the role of such enrichment functions, in this study, is to aid in representing the discontinuity surface independent of the mesh. Here, an enriched four noded $\mathcal{C}^o$ shear flexible quadrilateral plate element based on the consistency approach is used to analyse the flutter behaviour. 
The numerical integration over the elements intersected by the discontinuity surface is done based on conformal mapping~\citep{natarajanmahapatra2010}, which eliminates the need to sub-triangulate the elements intersected by the discontinuity surface. The influence of the crack length, the crack orientation, the flow direction and the material property on the evaluation of critical speed and the type of fluttering instability is numerically studied. 

The paper is organized as follows, the next section will give an introduction to FGM and a brief overview of Reissner-Mindlin plate theory. Section~\ref{numerics} presents results for the flutter analyses of cracked functionally graded material panels, followed by concluding remarks in the last section.

\section{Formulation}
\begin{figure}
\centering
\input{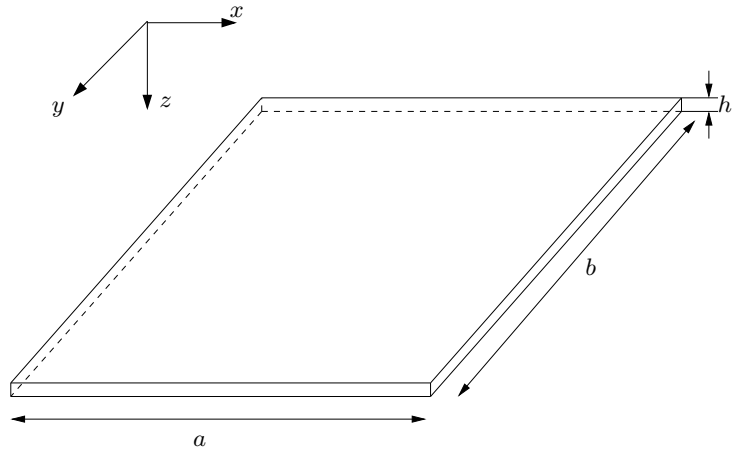}
\caption{Co-ordinate system of a rectangular FGM plate.}
\label{fig:platefig}
\end{figure}

A rectangular plate made of a mixture of ceramic and metal is considered with the coordinates $x,y$ along the in-plane directions and $z$ along the thickness direction (see \fref{fig:platefig}). The material on the top surface $(z=h/2)$ of the plate is ceramic and is graded to metal at the bottom surface of the plate $(z=-h/2)$ by a power law distribution. The effective Young's modulus $E$ and Poisson's ratio $\nu$ of the FGM evaluated using the rule of mixtures are:

\begin{eqnarray}
E = E_c V_c + E_m V_m \nonumber \\
\nu = \nu_c V_c + \nu_m V_m
\end{eqnarray}

where $V_i (i=c,m)$ is the volume fraction of the phase material. The subscripts $c$ and $m$ refer to ceramic and metal phases, respectively. The volume fraction of ceramic and metal phases are related by $V_c + V_m = 1$ and $V_c$ is expressed as:

\begin{equation}
V_c(z) = \left( \frac{2z+h}{2h} \right)^k 
\end{equation}

where $k$ is the volume fraction exponent $(k \geq 0)$, also called as the gradient index. The variation of the composition of ceramic and metal is linear for $k=$1, the value of $k=$ 0 represents a fully ceramic plate and any other value of $k$ yields a composite material with a smooth transition from ceramic to metal.

Using Mindlin formulation, the displacements $(u,v,w)$ at a point $(x,y,z)$ in the plate from the medium surface are expressed as functions of midplane displacements $(u_o,v_o,w_o)$ and independent rotations $\theta_x$ and $\theta_y$ of the normal in $xz$ and $yz$ planes, respectively, as:

\begin{eqnarray}
u(x,y,z,t) &=& u_o(x,y,t) + z \theta_x(x,y,t) \nonumber \\
v(x,y,z,t) &=& v_o(x,y,t) + z \theta_y(x,y,t) \nonumber \\
w(x,y,z,t) &=& w_o(x,y,t) 
\label{eqn:displacements}
\end{eqnarray}

The midplane membrane strains $\bveps_p$, bending strain $\bveps_b$ and shear strain $\bveps_s$ in are written as

\begin{eqnarray}
\renewcommand{\arraystretch}{1.5}
\bveps_p = \left\{ \begin{array}{c} u_{o,x} \\ v_{o,y} \\ u_{o,y}+v_{o,x} \end{array} \right\}, \hspace{0.25cm}
\renewcommand{\arraystretch}{1.5}
\bveps_b = \left\{ \begin{array}{c} \theta_{x,x} \\ \theta_{y,y} \\ \theta_{x,y}+\theta_{y,x} \end{array} \right\}, \hspace{0.25cm}
\renewcommand{\arraystretch}{1.5}
\bveps_s = \left\{ \begin{array}{c} \theta _x + w_{o,x} \\ \theta _y + w_{o,y} \end{array} \right\}. \hspace{1cm}
\renewcommand{\arraystretch}{1.5}
\end{eqnarray}

where the subscript `comma' represents the partial derivative with respect to the spatial coordinate succeeding it. The strain energy of the plate can be expressed in terms of the field variables $\boldsymbol{\delta} = (u_o,v_o,w_o,\theta_x,\theta_y)$ and their derivatives as:

\begin{equation}
U(\boldsymbol{\delta}) = {1 \over 2} \int_{\Omega} \left\{ \bveps_p^{\textup{T}} \mathbf{A} \bveps_p + \bveps_p^{\rm T} \mathbf{B} \bveps_b + 
\bveps_b^{\textup{T}} \mathbf{B} \bveps_p + \bveps_b^{\textup{T}} \mathbf{D}_b \bveps_b +  \bveps_s^{\textup{T}} \mathbf{E} \bveps_s\right\}~ \rmd \Omega
\label{eqn:potential}
\end{equation}

where the matrices $\mathbf{A}, \BB, \DD$ and $\mathbf{E}$ are the extensional, bending-extensional coupling, bending and transverse shear stiffness coefficients. The kinetic energy of the plate is given by:

\begin{equation}
T(\boldsymbol{\delta}) = {1 \over 2} \int_{\Omega} \left\{I_o (\dot{u}_o^2 + \dot{v}_o^2 + \dot{w}_o^2) + I_1(\dot{\theta}_x^2 + \dot{\theta}_y^2) \right\}~\rmd \Omega
\label{eqn:kinetic}
\end{equation}

where $I_o = \int_{-h/2}^{h/2} \rho(z)~dz, ~ I_1 = \int_{-h/2}^{h/2} z^2 \rho(z)~dz$ and $\rho(z)$ is the mass density that varies through the thickness of the plate. The work done by the applied non-conservative loads is:
\begin{equation}
W(\boldsymbol{\delta}) = \int_{\Omega} \Delta p w ~\rmd \Omega
\label{eqn:aerowork}
\end{equation}

where $\Delta p$ is the aerodynamic pressure. The aerodynamic pressure based on first-order, high Mach number approximation to linear potential flow is given by:

\begin{equation}
\Delta p = \frac{\rho_a U_a^2}{\sqrt{M_\infty^2 - 1}} \left[ \frac{\partial w}{\partial x} \cos \theta^\prime + \frac{\partial w}{\partial y} \sin \theta^\prime + \left( \frac{1}{U_a} \right) \frac{M_\infty^2 - 2}{M_\infty^2 - 1} \frac{\partial w}{\partial t} \right]
\label{eqn:aeropressure}
\end{equation}

where $\rho_a, U_a, M_\infty$ and $\theta^\prime$ are the free stream air density, velocity of air, Mach number and flow angle, respectively. The static aerodynamic approximation for Mach numbers between $\sqrt{2}$ and $2$ is~\citep{Dixon1966,Birman1990}:

\begin{equation}
\Delta p = \frac{\rho_a U_a^2}{\sqrt{M_\infty^2 - 1}} \left[ \frac{\partial w}{\partial x} \cos \theta^\prime + \frac{\partial w}{\partial y} \sin \theta^\prime  \right]
\label{eqn:aeropressurestat}
\end{equation}

The governing equations obtained using the minimization of total potential energy are solved based on the finite element method. The finite element equations thus derived are:

\begin{equation}
\left[ \left( \KK + \lambda \overline{\mathbf{A}}\right) - \omega^2 \bm\right] \boldsymbol{\delta} = \mathbf{0}
\label{eqn:finaldiscre}
\end{equation}

where $\KK$ is the stiffness matrix, $\bm$ is the consistent mass matrix, $\lambda = \frac{\rho_a U_a^2}{\sqrt{M_\infty^2 - 1}}$, $\overline{\mathbf{A}}$ is the aerodynamic force matrix and $\omega$ is the natural frequency. When $\lambda = $0, the eigenvalue of $\omega$ is real and positive, since the stiffness matrix and mass matrix are symmetric and positive definite. However, the aerodynamic matrix $\overline{\mathbf{A}}$ is unsymmetric and hence complex eigenvalues $\omega$ are expected for $\lambda >$ 0. As $\lambda$ increases monotonically from zero, two of these eigenvalues will approach each other and become complex conjugates. In this study, $\lambda_{cr}$ is considered to be the value of $\lambda$ at which the first coalescence occurs.

\section{Results and discussion}
\label{numerics}

In this section, we present the critical aerodynamic pressure and the critical frequency of a cracked simply supported FGM panels using the extended Q4 formulation~\citep{natarajanbaiz2011}. The element has five nodal degrees of freedom $(u_o,v_o,w_o,\theta_x,\theta_y)$. The formulation includes transverse shear deformation and in-plane and rotary inertia effects. A full integration scheme is applied to evaluate the various strain energy terms~\citep{Somashekar1987,Ganapathi1991,natarajanbaiz2011}. A simply supported boundary condition is assumed for the current study, given by:

\begin{eqnarray}
u_o = w_o = \theta_y = 0 \hspace{1cm} ~\textup{on} ~ x=0,a \nonumber \\
v_o = w_o = \theta_x = 0 \hspace{1cm} ~\textup{on} ~ y=0,b
\end{eqnarray}

 In all cases, we present the non dimensionalized critical aerodynamic pressure, $\lambda_{cr}$ and critical frequency $\omega_{cr}$ as, unless specified otherwise:

\begin{eqnarray}
\Omega_{cr} = \omega_{cr} a^2 \sqrt{ \frac{\rho_c h}{D_c}} \nonumber \\
\lambda_{cr} = \lambda_{cr} \frac{a^3}{D_c}
\label{eqn:nondimfreq}
\end{eqnarray}

where $D_c = {E_c h^3 \over 12(1-\nu_c^2)}$ is the bending rigidity of the plate, $E_c, \nu_c$ are the Young's modulus and Poisson's ratio of the ceramic material and $\rho_c$ is the mass density. For this study, the plate thickness is assumed to be $a/h=$ 100, unless specified otherwise. In order to be consistent with the existing literature, properties of the ceramic are used for normalization. The effect of material property, the crack orientation $\theta$, the crack length $d/a$ and the flow angle $\theta^\prime$ on the flutter behaviour are studied. Based on progressive mesh refinement, a $34\times34$ structured mesh is found to be adequate to model the full plate for the present analysis. The material properties used for the FGM components are listed in Table~\ref{table:tempdepprop}. 

Before proceeding with the detailed study, the formulation developed herein is validated against available results pertaining to the critical aerodynamic pressure and critical frequency for an isotropic plate with and without a crack. The computed critical aerodynamic pressure and the critical frequency for an isotropic square plate with various boundary conditions is given in  Table~\ref{table:isoFlutComp}. Table~\ref{table:isoFlutCompCrack} gives a comparison of the computed frequency and the aerodynamic pressure for a cantilevered square plate with a side crack immersed in a normal flow $(\theta^\prime=$ 0$^\circ)$. It can be seen that the numerical results from the present formulation are found to be in  good agreement with the existing solutions.  \\

\begin{figure}[htpb]
\centering
\input{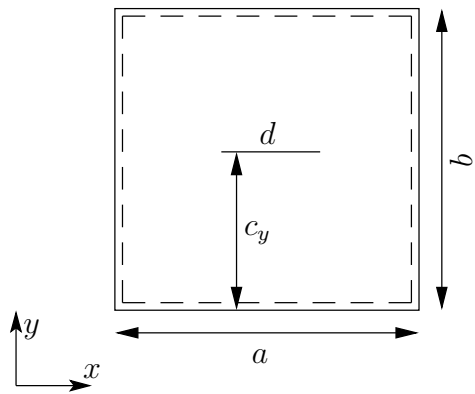}
\caption{Simply supported plate with a center crack.}
\label{fig:SS2Plate}
\end{figure}

\begin{figure}[htpb]
\centering
\includegraphics[scale=0.8]{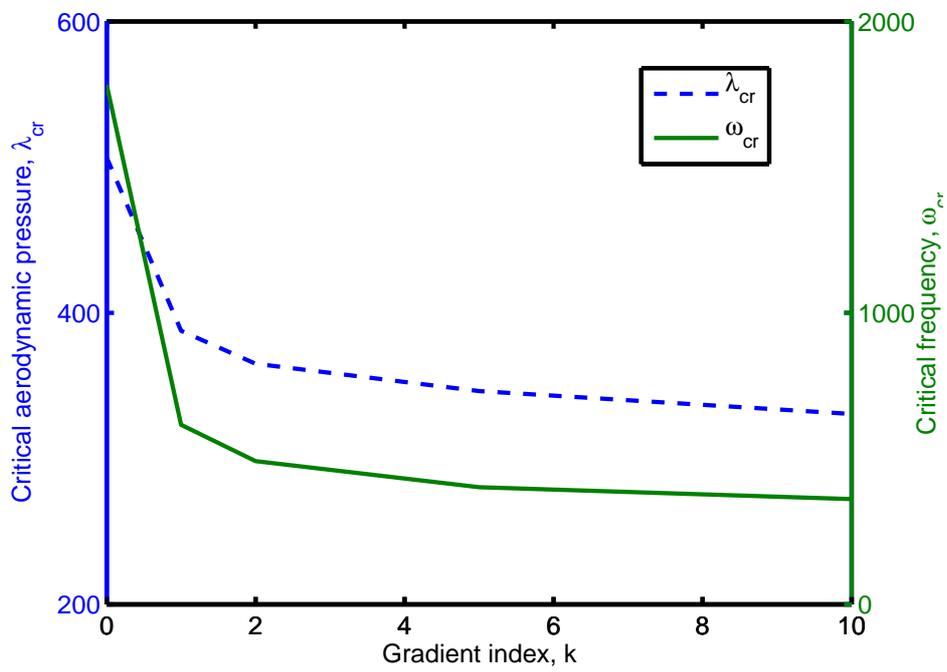}
\caption{Variation of critical aerodynamic pressure $\lambda_{cr}$ and the critical frequency $\omega_{cr}$ with gradient index in a normal flow $\theta^\prime =$0$^\circ$ for a simply supported square FGM plate with a center horizontal crack with $d/a=$ 0.5 }
\label{fig:pressvsgindex}
\end{figure}

\begin{figure}[htpb]
\centering
\includegraphics[scale=0.8]{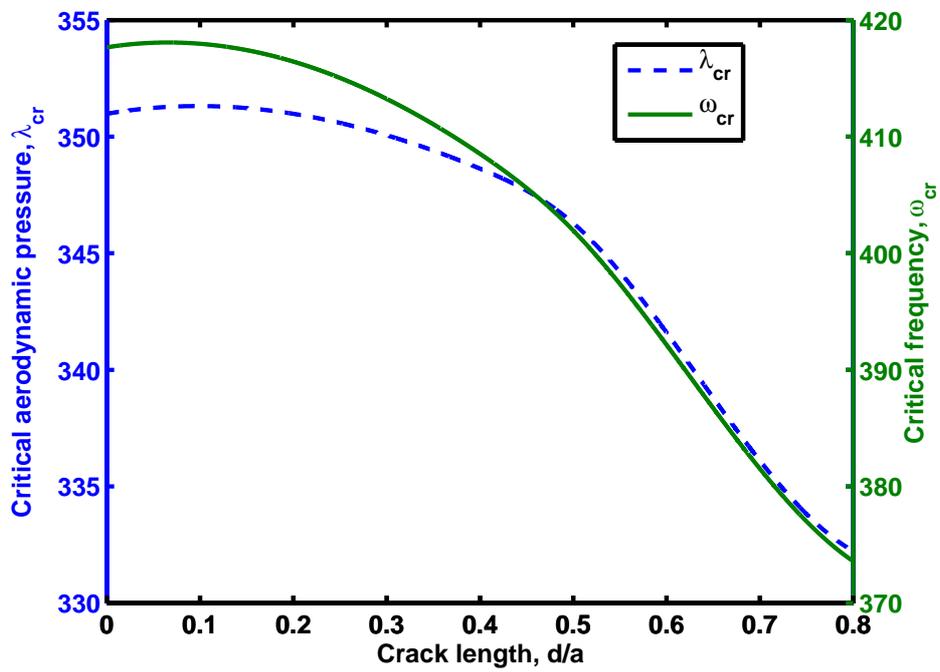}
\caption{The critical aerodynamic pressure $\lambda_{cr}$ and the critical frequency $\omega_{cr}$ as a function of crack length $d/a$ for a simply supported square FGM plate with a center horizontal crack in a normal flow $\theta^\prime=$0 $^\circ$ with gradient index $k=$ 5.}
\label{fig:pressvccrklen}
\end{figure}

Next, the flutter characteristics of square simply supported cracked FGM plates made of silicon nitride (Si$_3$N$_4$) and steel (SUS304) is investigated. Consider a plate of uniform thickness, $h$ and with length and width as $a$ and $b$, respectively. \fref{fig:SS2Plate} shows a plate with all edges simply supported with a center crack of length $d$ and a distance of $c_y$ from the $x-$ axis. In this example, the influence of the crack length $d/a$, the crack orientation $\theta$, the flow angle $\theta^\prime$ and the gradient index $k$ on the critical aerodynamic pressure and the critical frequency is studied. \fref{fig:pressvsgindex} shows the variation of critical aerodynamic pressure $\lambda_{cr}$ and critical frequency $\omega_{cr}$ with gradient index for a center horizontal crack of length $d/a=$ 0.5 in a normal flow $\theta^\prime=$ 0$^\circ$. It can be seen that with increasing gradient index, both the critical aerodynamic pressure and the critical frequency decreases. This is because of the stiffness degradation due to increase in metallic volume fraction. \fref{fig:pressvccrklen} shows the influence of the crack length $d/a$ on the critical aerodynamic pressure and the critical pressure. It is observed that as the crack length increases, the critical aerodynamic pressure and the critical frequency decreases. This is due to the fact that increasing the crack length increases the local flexibility and thus decreases the frequency. It can be seen that the combined effect of increasing the crack length and the gradient index is to lower the critical frequency and the critical pressure. It can also be seen from \fref{fig:pressvsgindex} that the critical pressure and the frequency drops sharply for small increase in the metallic volume fraction, but with further increase in the metallic volume fraction, the drop in the pressure and the frequency is marginal. In both the cases, the decrease in critical aerodynamic pressure and the critical frequency is due to the stiffness degradation.

\begin{figure}[htpb]
\centering
\includegraphics[scale=0.8]{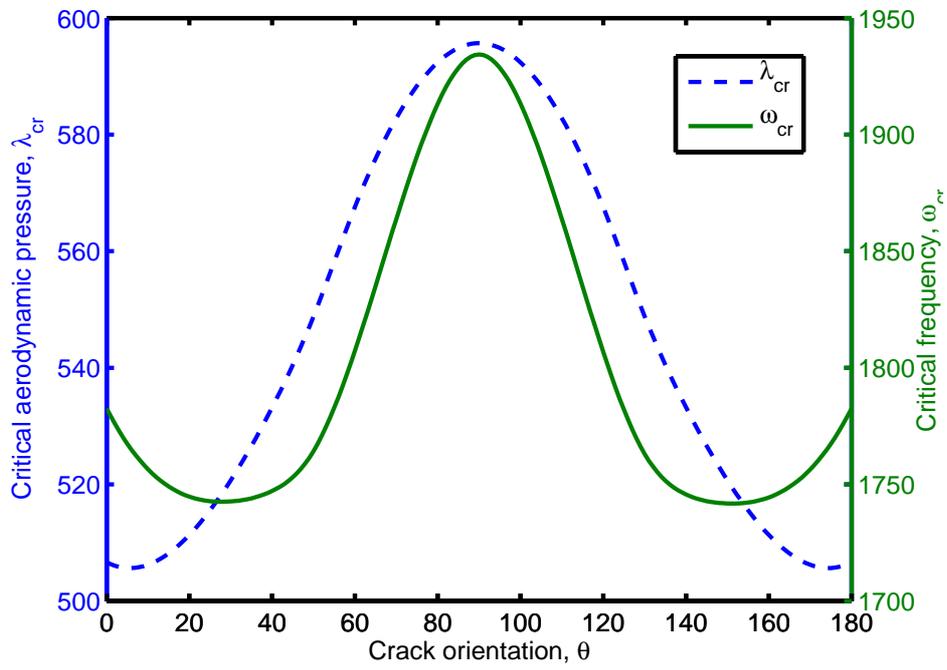}
\caption{Critical aerodynamic pressure $\lambda_{cr}$ as a function of crack orientation for a simply supported square FGM plate with a center crack $d/a=$ 0.5 in a normal flow $\theta^\prime=$ 0 and gradient index $k=$ 0.}
\label{fig:crkanglevslambda}
\end{figure}

\fref{fig:crkanglevslambda} shows the influence of crack angle $\theta$ on the critical pressure and the frequency for a simply supported square FGM plate with gradient index $k=$ 0, immersed in a normal flow $\theta^\prime=$0$^\circ$. From \fref{fig:crkanglevslambda} it can be seen that with increase in the crack orientation, the critical pressure increases gradually until the crack is oriented at right angles to the flow angle and with further increase in the crack orientation, the critical pressure and frequency decreases. Further, it is observed that the critical pressure is lowest for a crack orientation $\theta = $ 0$^\circ$ and 180$^\circ$. At these crack orientations, the crack is aligned to the flow direction. The critical frequency and the pressure values tend to be symmetric with respect to a crack orientation $\theta=$90$^\circ$.

\begin{figure}[htpb]
\centering
\includegraphics[scale=0.8]{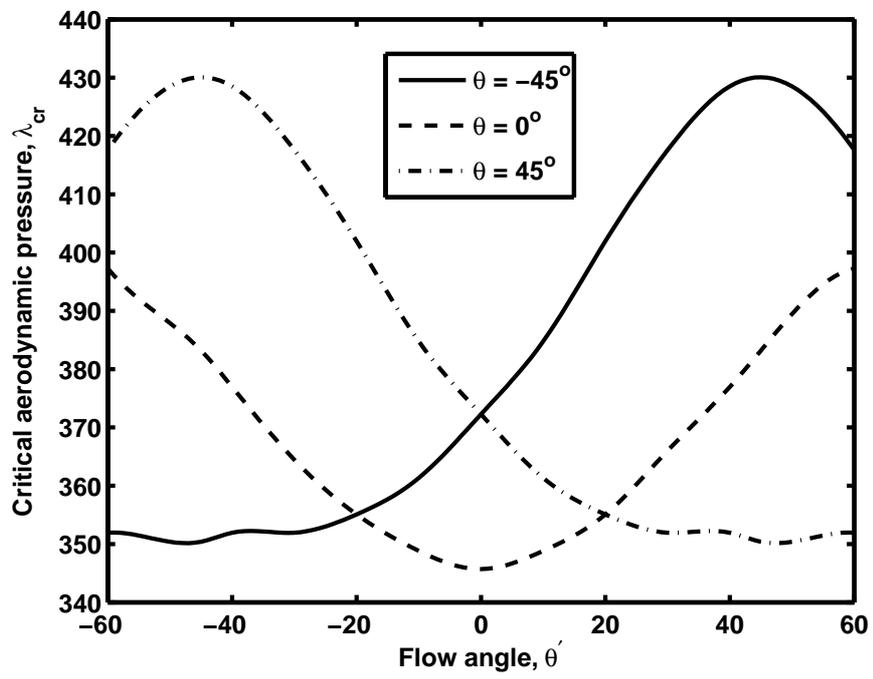}
\caption{Effect of flow angle $\theta^\prime$ on the critical aerodynamic pressure for a simply supported square FGM plate with a center horizontal crack with $d/a=$ 0.5 and gradient index $k=$ 5.}
\label{fig:fgmFlowAngle}
\end{figure}

The influence of the flow angle on the critical pressure for a simply supported square FGM plate with gradient index $k=$ 5 is shown in \fref{fig:fgmFlowAngle} for different crack orientations, i.e., $\theta = $ 0$^\circ$, $\pm$45$^\circ$. It is observed that the critical pressure is minimum when the crack is aligned to the flow angle. The critical frequency increases for other flow angles irrespective of the orientation of the crack.

\section{Conclusion}
The flutter characteristics of cracked FGM panels immersed in a supersonic flow has been analyzed based on first-order shear deformation theory through the extended finite element approach. The aerodynamic force is accounted for assuming the first-order Mach number approximation potential flow theory and the homogenized material properties are estimated by rule of mixtures. Numerical experiments have been conducted to bring out the effect of the gradient index, the crack length, the crack orientation and the flow angle on the flutter characteristics of the FGM plate. From the numerical study, it can be concluded that with increase in both the gradient index and the crack length, the critical pressure and the frequency decreases. In both cases, the decrease is due to the stiffness degradation. It is also observed that the critical frequency and the pressure is minimum when the crack is aligned to the flow.

\bibliographystyle{abbrvnat}
\bibliography{myFlut}

\newpage

\begin{figure}[htpb]
\centering
\input{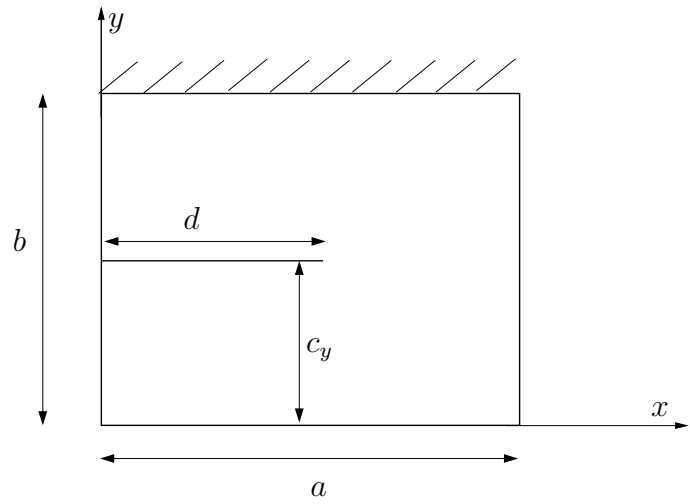}
\caption{Cantilevered plate with a side crack: geometry}
\label{fig:cantiplate}
\end{figure}

\begin{table}
\renewcommand\arraystretch{1.5}
\caption{Temperature dependent coefficient for material Si$_3$N$_4$/SUS304, Ref~\citep{Reddy1998,Sundararajan2005}.}
\centering
\begin{tabular}{lcccccc}
\hline
Material & Property & $P_o$ & $P_{-1}$ & $P_1$ & $P_2$ & $P_3$  \\
\hline
\multirow{2}{*}{Si$_3$N$_4$} & $E$(Pa) & 348.43e$^9$ &0.0& -3.070e$^{-4}$ & 2.160e$^{-7}$ & -8.946$e^{-11}$  \\
& $\alpha$ (1/K) & 5.8723e$^{-6}$ & 0.0 & 9.095e$^{-4}$ & 0.0 & 0.0 \\
\cline{2-7}
\multirow{2}{*}{SUS304} & $E$(Pa) & 201.04e$^9$ &0.0& 3.079e$^{-4}$ & -6.534e$^{-7}$ & 0.0  \\
& $\alpha$ (1/K) & 12.330e$^{-6}$ & 0.0 & 8.086e$^{-4}$ & 0.0 & 0.0 \\
\hline
\end{tabular}
\label{table:tempdepprop}
\end{table}

\begin{table}[htpb]
\renewcommand\arraystretch{1.5}
\caption{Comparison of critical aerodynamic pressure and coalescence frequency for an isotropic plate with various boundary conditions $(a/b=1, a/h=100, \nu=0.3, \theta^\prime = 0)$.}
\centering
\begin{tabular}{cccc}
\hline 
Reference & Flutter bounds & \multicolumn{2}{c}{Boundary condition} \\
\cline{3-4}
1& 1 & Simply supported & Clamped \\
\hline
\multirow{2}{*}{Ref.~\citep{Sander1973}} & $\lambda_{cr}$ & 512.2 & 853.40 \\
& $\omega_{cr}$ & 1844.00 & 4292.00 \\
\multirow{2}{*}{Ref.~\citep{Han1983}} & $\lambda_{cr}$ & 512.33 & 852.73\\
& $\omega_{cr}$ & 1840.55 & 4294.07 \\
\multirow{2}{*}{Ref.~\citep{Lin1989}} & $\lambda_{cr}$ & 512.58 & 851.50\\
& $\omega_{cr}$ & 1847.50 & 4290.00 \\
\multirow{2}{*}{Ref.~\citep{Prakash2006}} & $\lambda_{cr}$ & 511.11 & 852.34 \\
& $\omega_{cr}$ & 1840.29 & 4274.32 \\
\multirow{2}{*}{Present} & $\lambda_{cr}$ &513.48  & 854.80\\
& $\omega_{cr}$ & 1849.50 &4297.00\\
\hline
\end{tabular}
\label{table:isoFlutComp}
\end{table}

\begin{table}[htpb]
\renewcommand\arraystretch{1.5}
\caption{Comparison of flutter solutions for a square isotropic cantilevered panel with an edge crack $d/a = 0.3$ (see \fref{fig:cantiplate}) in a normal flow $\theta^\prime = 0$ with $a/h=100, \nu=0.3$. }
\centering
\begin{tabular}{ccc}
\hline 
Solutions & \multicolumn{2}{c}{Flutter bounds} \\
\cline{2-3}
& $\omega_{cr}$ & $\lambda_{cr}$ \\
\hline 
Ref.~\citep{Lin1989} & 22.60 & 46.30 \\
Ref.~\citep{Lin1991} & 21.43 & 46.50 \\
Present & 22.35 & 45.50 \\
\hline
\end{tabular}
\label{table:isoFlutCompCrack}
\end{table}

\end{document}